\documentclass[12pt]{amsart}
\usepackage{epsfig,color}

\newtheorem{Thm}[equation]{Theorem}
\newtheorem{Cor}[equation]{Corollary}
\newtheorem{Lem}[equation]{Lemma}
\newtheorem{Pro}[equation]{Proposition}

\theoremstyle{definition}

\newtheorem{Def}[equation]{Definition}
\newtheorem{Exa}[equation]{Example}

\theoremstyle{remark}
\newtheorem{Que}[equation]{\bf Question}
\newtheorem{Rem}[equation]{Remark}

\numberwithin{equation}{section}
\makeatletter
\renewcommand{\c@figure}{\c@equation}
\makeatother

\newcommand{\maths}[1]{{\bf #1}}

\newcommand{\RR}{\maths{R}}
\newcommand{\NN}{\maths{N}}

\renewcommand{\SS}{\maths{S}}

\newcommand{\ra}{\rightarrow}

\newcommand{\bord}{\partial}

\newcommand{\dem}{{\sl {Proof.} \  }}
\newcommand{\eop}[1]{{\flushright\hfill\fbox{\bf #1}}}

\newcommand{\ack}{\noindent{\bf Acknowledgement.}}

\begin{document}

\title[A closed topological cube]{Applications of Three Dimensional Extremal Length, 
 I: tiling of a topological cube}

\author{Sa'ar Hersonsky}

\address{Department of Mathematics\\ 
University of Georgia\\ 
Athens, GA 30602}

\urladdr{http://www.math.uga.edu/~saarh}
\email{saarh@math.uga.edu}
\date{}

\begin{abstract}
 Let $\mathcal{T}$ be a triangulation of a closed topological cube $Q$, and let $V$ be the set of vertices of $\mathcal{T}$. Further assume that the triangulation satisfies a technical condition which we call the {\it triple intersection property} (see Definition~\ref{de:intersecting}). Then there is an essentially unique tiling ${\mathcal C}=\{C_v  : v\in V\}$ of a rectangular parallelepiped $R$ by cubes, such that for every edge $(u,v)$ of  $\mathcal{T}$ the corresponding cubes $C_v,C_u$ have nonempty intersection, and such that the vertices corresponding to the cubes at the corners of $R$ are at the corners of $Q$.  
Moreover,  the sizes of the cubes are obtained as a solution of a variational problem which is a discrete version of the notion of extremal length in $\RR^3$. 
 \end{abstract}

\maketitle

\section{Introduction}
\label{se:Intro}

One exciting connection between two dimensional conformal geometry and packing is provided by the circle packing theorem, which asserts that if $G$ is a finite planar graph, then there is a packing of Euclidean disks in the complex plane whose contact graph is $G$. This remarkable result was first proved by Koebe (\cite{Ko}) as a consequence of his theorem that every finite planar domain is conformally equivalent to a circle domain. Koebe's result was rediscovered and vastly generalized by Thurston (\cite{Th1},\cite[Chapter 13]{Th2}) as a corollary of Andreev's Theorem (\cite{An1},\cite{An2}).  Thurston also conjectured that a sequence of maps, naturally associated to circle packings of a simply connected domain, converges to the Riemann map from this domain to the unit disk. This conjecture was proved by Rodin and Sullivan (\cite{RoSu}), thus providing a second foundational connection between 
circle packing and conformal maps. In his thesis and later work, Schramm (see in particular \cite[Theorem 6.1]{Sch2}) generalized the circle packing theorem allowing the tiles in the packing to be homothetic to $C^1$ closed topological disks. His result is based on an elaborated conformal uniformization theorem of Brandt and Harrington, and established another important connection between circle packing and two dimensional conformal geometry. 

Following a suggestion by Thurston, Schramm  (\cite{Sch1}) studied the case in which the tiles in the packing are squares. An independent and similar study was carried out by Cannon, Floyd and Parry (\cite{CaFlPa}), as part of their attempts to resolve Cannon's conjecture. Both results are based on  {\it discrete extremal length} arguments, a notion first developed by Cannon (\cite{Ca}). This notion has its origin in the subject of two dimensional quasiconformal maps, where extremal length arguments are essential. There are a wealth of other results, relating combinatorics and packing that one should mention. Benjamini and Schramm (\cite{BeSch}) studied the case where the tiled set is an infinite straight cylinder while Kenyon (\cite{Ke}) allowed the tiles of the packing to be polygons.  It is important to recall that the results in \cite{Ke}, as well as more contemporary work by the author of this paper (\cite{Her1},\cite{Her2}), do not use discrete extremal length methods. Rather, the usage of {\it discrete harmonic functions} as  
first employed by Dehn (\cite{De}), and later on by Brooks, Smith, Stone and Tutte (\cite{BSST}), is utilized. These results are also different from the ones obtained in \cite{Sch1}, \cite{CaFlPa} and \cite{BeSch}. In these papers, a tile corresponds to a vertex in a given triangulation. In \cite{De}, \cite{BSST}, \cite{Ke}, \cite{Her1} and \cite{Her2}, a tile corresponds to an edge of the triangulation. It is also worth noting that in \cite{Her1} and \cite{Her2}, multi-connected, bounded,  planar domains were studied (under the framework of boundary value problems on graphs) for the first time.  This provides tiling of higher genus surfaces with conical singularities by rectangles.

\medskip

The main goal of this paper is to provide the first  connection between discrete extremal length in $\RR^{3}$, a notion which we will recall in \S\ref{se:section1}, and tiling by cubes (see  Remark~\ref{re:gene} for one possible generalization). In this paper, we study the case of a topological closed cube. It is interesting to note that quite recently Benjamini and Schramm (\cite{BeSch1}),  as well as Benjamini and Curien (\cite{BeCu}) have explored different and interesting applications of discrete extremal length in higher dimensions. 

\medskip

Before stating the main result of this paper we make 
\begin{Def} 
\label{de:trianofcube}
Let $B$ be a closed triangulated topological ball, and let $V, E$ and $F$ denote the set of vertices, edges, and faces of the triangulation, respectively. Let $\partial B = B_1\cup \bar{B_1}\cup B_2\cup \bar{B_2}\cup B_3\cup \bar{B_3}$ be a decomposition of $\partial B$ in such a way that each $B_j$ is a nonempty connected union of faces of the triangulation, $B_i\cap \bar{B_j}=\emptyset$ for $i=j$, and consists of a union of edges of the triangulation, if $i\neq j$. The collection ${\mathcal T}=\{V,E,F;B_1,\bar{B_1},B_2,\bar{B_2},B_3,\bar{B_3}\}$ will be called a triangulation of a topological cube.
We will denote by $B_1$ the \it{base face}, by $\bar {B_1}$ the \it{top face}, by $B_2$ the \it{front face}, by $\bar {B_2}$ the {\it back face}, by $B_3$ the {\it left face}, and by $\bar{B_3}$ the \it{right face}. 
\end{Def}

The main result of this paper is 
\begin{Thm}
\label{th:main}
Let ${\mathcal T}=\{V,E,F;B_1,\bar{B_1},B_2,\bar{B_2},B_3,\bar{B_3}\}$ be a triangulation of a topological cube which has the triple intersection property. Then there exists a positive number $h$ and a cube tiling  ${\mathcal C}=\{C_v : v\in V\}$ of the rectangular parallelepiped $R=[0,\sqrt{h^{-1}}]\times[0,\sqrt{h^{-1}}]\times[0,h]$ such that 
\begin{equation}
\label{eq:comb}
C_v\cap C_u\neq\emptyset \ \mbox{\rm whenever}\ (v,u)\in E.
\end{equation}

In addition, let $R_1,\bar{R_1},R_2,\bar{R_2}, R_3$, and $\bar{R_3}$ be the base, top, front, back, left, and right faces of $R$, respectively.  Then it can also be arranged that for $i=1,2,3$ we have
\begin{equation}
\label{eq:bound}
C_v\cap R_i(\bar R_i)\neq\emptyset \ \mbox{\rm whenever}\ v\in B_i(\bar {B_i}).
\end{equation}
Under these conditions, the number $h$ and the tiling ${\mathcal C}$ are uniquely determined.
\end{Thm}

Theorem~\ref{th:main} is a generalization to three dimensions of the main result (Theorem 1.3) as well as the techniques in \cite{Sch1}, under an extra assumption. 
Schramm's proof (along with the proof given by Cannon, Floyd and Parry) fails to work in three dimension. The planarity of the triangulation is essential in their proofs.
The triple intersection property (Definition~\ref{de:intersecting}), enables us to extract the ideas and techniques in \cite{Sch1} and carry out our proof, which then becomes straightforward.
 We have examples in which this property holds, and  which can be characterized by saying that the  triangulation has a {\it spine} 
(Definition~\ref{de:narrowneck}).

The rest of this paper is organized as follows. In \S\ref{se:section1}, we recall the notion of discrete extremal length in dimension three and prove the existence and uniqueness of an {\it extremal metric} (up to scaling). In \S\ref{se:section2}, we prove that tiling by cubes induces an extremal metric, and in \S\ref{se:section3}, we prove that an extremal metric induces a tiling by cubes. Finally, \S\ref{se:section4} is devoted to  questions and suggestions for further research.

\medskip

\ack {\small \  After receiving an early version of this paper, I.~Benjamini brought to our attention recent references regarding applications of discrete extremal length in high dimensions (\cite{BeSch1},\cite{BeCu}). 

\section{Perspective and basic definitions}
\label{se:section1}
In the 1940's, Ahlfors and Beurling have refined existence methods (by Gr\"{o}tzch and Teichm\"{u}ller) and extremal length was used as a conformally invariant measure of planar curve families (see for instance \cite{Ah} for a useful account). 
 L\"owener (\cite{Lo}) showed how this method can be extended to three dimensions by defining a conformal capacity for rings in Euclidean $3$-space by means of a Dirichlet integral.  V\"{a}is\"{a}l\"{a} 
(\cite{Va1,Va2}) and \u{S}abat (\cite{Sa1,Sa2}) have used extremal length arguments to study quasiconformal mappings in $3$-space, each of them has introduced a new kind of capacity for a ring in three dimensional Euclidean space. Shortly afterwards, Gehring (\cite[Theorem 1]{Ge}) showed that slight modifications of their definitions are equivalent to the one given by L\"owener. Extremal length arguments have continued to be useful tools in the theory of quasiconformal mappings of the plane and have found profound application in Teichm\"{u}ller theory and the theory of hyperbolic manifolds and their deformations.
  We wish to restrict the background and preliminaries to a minimum. Hence, we will not describe many of the modern definitions and exciting applications of extremal length in the general setting of metric measure spaces. The interested reader is advised to consult for example Heinonen \cite{Hei}, and the references therein for an enjoyable and extensive account.

Our main result (Theorem~\ref{th:main}) generalizes the main results of \cite{Sch1} and \cite{CaFlPa} that are based on 
 Cannon's definition of two dimensional extremal length on a graph (see \cite{Ca}).  Cannon's definition was extended for arbitrary graphs in \cite[Section 9]{Sch1}.  The definition below is a special case of the one given in \cite[Section 9]{Sch1}, and is suitable to the applications of this paper. It is essentially a discrete version of the definition given by V\"{a}is\"{a}l\"{a}.  Let $G=(V,E)$ be a finite connected graph.  A {\it path} in $G$ is a sequence of vertices $\{v_0,v_1,\ldots,v_k\}$ such that any two successive vertices are connected by an edge. A nonnegative function $m:V\rightarrow [0,\infty)$ will be called a {\it metric} on $G$. Given a path $\alpha= \{v_0,v_1,\ldots,v_k\}$ in G and a metric $m$, we define the $m$-length of $\alpha$ 
as
\begin{equation}
\label{eq:length}
l_m(\alpha)= \sum_{i=0}^{k}m(v_i).
\end{equation}
Given $A_1,A_2\subset V$, we define their $m$-distance to be
\begin{equation}
\label{eq:distsubsets}
l_m(A_1,A_2)=\inf_{\alpha}l_m(\alpha),
\end{equation}
where the infimum is taken over all paths $\alpha$ which join $A_1$ to $A_2$. The {\it volume} of the metric $m$ is defined to be the cube of its 3-norm, i.e
\begin{equation}
\label{eq:vol}
\mbox{\rm vol}(m)= ||m||_3^3=\sum_{v\in V} m(v)^3, 
\end{equation}   
and the {\it normalized length} of $(m,A_1,A_2)$ is defined as
\begin{equation}
\label{eq:norlength}
{\hat l_m}=\frac{l_m^3}{\mbox{\rm vol}(m)}.
\end{equation}
Finally, the {\it extremal three dimensional length} of $(G,A_1,A_2)$ is defined  by 
\begin{equation}
\label{eq:exlength}
\lambda(G; A_1,A_2)= \sup_{m} {\hat l_m},
\end{equation}
where the supremum is taken over the set of all metrics $m$ with positive volume. An {\it extremal metric} for $(G;A_1,A_2)$ is one which realizes this supremum. Observe that for any positive constant $c$, $l_{cm}=c\, l_{m}$. Hence, ${\hat l_m}$ is a {\it conformal invariant} in this discrete setting.

\medskip

An important theorem which asserts the existence and uniqueness (up to scaling) of  an extremal metric was proved independently by Schramm and by Cannon, Floyd and Parry. We now recall the proof given by Cannon, Floyd and Parry (\cite[Theorem 2.2.1]{CaFlPa}), which applies (with negligible modifications) to our setting, in order to make this paper self-contained.   

\begin{Thm} 
\label{th:existence}
There is a unique extremal metric $m_0$ for $(G, A_1, A_2)$ such that $\mbox{\rm vol}(m_0)=1$.
\end{Thm}
\begin{proof}
Let $\NN$ denote the set of natural numbers, and let $P$ denote a nonempty finite subset of $\NN^n\setminus \{0\}$, where $n$ is the cardinality of $V$.  A path  in $G$ corresponds to an element in $P$. A metric $m$ on $G$ corresponds to a vector (which we will keep denoting by $m$) $m=(m_1,m_2,\ldots, m_n)\in \RR^n\setminus \{0\}$, with $m_i\geq 0$ for $i=1,\ldots n$. The length of a path with respect to the metric $m$ is then given by the standard scalar product in $\RR^n$. Since scaling does not change the extremal length $\lambda(G; A_1,A_2)$, we may restrict $m$ to the subset of vectors of $\SS^{n-1}$ in which each coordinate is nonnegative. Existence of an extremal metric now follows from the fact that the function which maps a metric $m\in \SS^{n-1}$ to $l_m(A_1,A_2)$ is the minimum of a finite number of continuous functions; hence, it is continuous. It now easily follows that $\lambda(G; A_1,A_2)$ is attained and is positive. Uniqueness essentially follows since 
$(\SS^{n-1},\| . \|_{3})$ is strictly convex. Given $m_1, m_2$ distinct nonnegative metrics in $\SS^{n-1}$ such that 
\begin{equation}
\label{eq:restriction}
l_{m_1}(A_1,A_2)\geq l_{m_2}(A_1,A_2)\  \mbox{\rm and}\  t\in (0,1),\ \mbox{\rm we let} 
\end{equation}

\begin{equation}
\label{eq:ontheline}
v= t m_1 + (1-t) m_2 \ \mbox{\rm(note that $0<||v||<1$)}.
\end{equation}
 Then for any path $p\in P$ we have

\begin{equation}
\begin{split}
       \frac{1}{||v||}<v,p>= & \frac{1}{||v||} (t<m_1,p>+ (1-t)<m_2,p>) \\
          \geq&\frac{1}{||v||}(t l_{m_1}(A_1,A_2)+ (1-t)l_{m_2}(A_1,A_2))  \\
          \geq&\frac{1}{||v||}l_{m_2}(A_1,A_2)>l_{m_2}(A_1,A_2),
    \end{split}
\end{equation}
which clearly shows that the extremal metric (up to scaling) is unique.

\end{proof}

\begin{Rem}
\label{re:optimalweight}
The geometry of discrete two dimensional extremal metrics was studied extensively by Parry and later on by Cannon, Floyd and Parry (see for example \cite{Ca},\cite{CaFlPa}).  We leave the study of the geometry of three dimensional extremal metrics for the future, since for the purposes of this paper only the assertion of Theorem~\ref{th:existence} is needed.
\end{Rem}

\begin{Rem}
An interesting recent reference by Wood (\cite{Woo}) explores some of the complications arising by the two inequivalent ways of carrying the notion of conformal modulus of a ring domain to a triangulated annulus. The first is by assigning a metric (as we do) to the vertices, and the second assigned a metric to the edges. It would be interesting to investigate this point in the three dimensional case.
\end{Rem}

\section{Cube tilings give extremal metrics}
\label{se:section2}
In this section, we prove that a cube tiling of a rectangular parallelepiped yields in a natural way an extremal metric. Our proof is carried out by modifying the main idea of the proof of Lemma 4.1 in \cite{Sch1} to three dimensions. In order to ease the notation, and since we are working with  fixed data, we let $l_x$ denote $l_{x}(B_1,\bar B_1)$ in the lemma below.
\begin{Lem}
\label{le:cutoex}
Let ${\mathcal T}=\{V,E,F;B_1,\bar{B_1},B_2,\bar{B_2},B_3,\bar{B_3}\}$ be as in Theorem~\ref{th:main}, and suppose that $h$ and ${\mathcal C}$ satisfy all the conditions there. Let $G=(V,E)$ be the 1-skeleton of ${\mathcal T}$, and let $s(v)$ denote the edge length of the square $C_v$. Then 
$s$ is an extremal metric for $(G, B_1,\bar{B_1})$.
\end{Lem}
\begin{proof}
Let $m$ be an arbitrary metric on $G$ with positive volume.
For every
\begin{equation}
\label{eq:domainof}
 (t,s)\in [0,\sqrt{h^{-1}}]\times[0,\sqrt{h^{-1}}],
  \end{equation}  
let
\begin{equation}
\label{eq:inverseofvertical}
\gamma_{t,s}=\left\{v\in V: \beta_{t,s}\cap C_v\neq\emptyset \right\}, \mbox{\rm where}\ \beta_{t,s}=(t,s)\times\RR. 
\end{equation}

\begin{figure}[htbp]
\begin{center}
 \scalebox{.44}{ \input{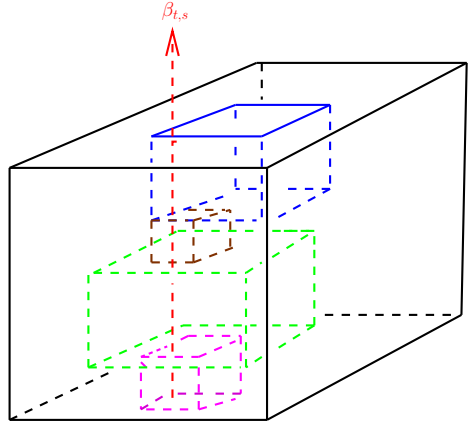}}
 \caption{$\beta_{t,s}$ going through cubes in the tiling.}
\label{figure:quad}
\end{center}
\end{figure}

It is clear that $\gamma_{t,s}$ contains a simple path in $G$ joining $B_1$ to $\bar{B_1}$. Hence, for every $(t,s)\in [0,\sqrt{h^{-1}}]\times[0,\sqrt{h^{-1}}]$ we have 

\begin{equation}
\label{eq:simplepath}
l_m\leq \sum_{v\in \gamma_{t,s}}m(v).
\end{equation}
We now integrate this inequality over $[0,\sqrt{h^{-1}}]\times[0,\sqrt{h^{-1}}]$ to obtain

\begin{equation}
\label{eq:integrate}
(\sqrt{h^{-1}})^{2}l_m \leq \int_{0}^{ \sqrt{h^{-1}} }\int_{0}^{\sqrt{h^{-1}}}\sum_{v\in \gamma_{t,s}}m(v)\,dt\,ds.
\end{equation}
Since ${\mathcal C}$ is a tiling and every $v\in V$ contributes $m(v)$ to the integral on the right hand side for an area measure of size $s(v)^{2}$, the integral can be rearranged to yield the following inequality
\begin{equation}
\label{eq:integrate1}
(\sqrt{h^{-1}})^{2}l_m \leq \sum_{v\in V}m(v)s(v)^2.
\end{equation}
We now let $s_1(v)=s(v)^2$ for every $v\in V$, and apply the Holder inequality with
$p=3$ and $q=3/2$ to the above inequality to obtain
\begin{equation}
\label{eq:holder}
l_m \leq h \sum_{v\in V}m(v)s(v)^2\leq h\, ||m||_{3}||s_1||_{3/2}= h\, (\sum_{v\in V}m(v)^{3})^{1/3} (\sum_{v\in V}s_1(v)^{3/2})^{2/3}.
 \end{equation}
Since $(\sum_{v\in V}s_1(v)^{3/2})^{2/3}= (\sum_{v\in V}s(v)^3)^{2/3}$ and $||s||_{3}^{3}=\mbox{\rm vol}(R)=1$, we finally obtain that
\begin{equation}
\label{eq:holder1}
l_m \leq h\, ||m||_{3}.
 \end{equation}

It is clear that $l_s=h$ and therefore that
\begin{equation}
\label{eq:endle}
{\hat l_m}=\frac{l_m^3}{||m||_3^3}\leq \frac{l_s^3}{||s||_3^3}={\hat l_s},
\end{equation} 
which implies that $s$ is extremal.

\end{proof}

\section{Extremal metrics give cube tiling}
\label{se:section3}

The main result of this section is Theorem~\ref{th:exttocubes} which asserts that an extremal metric, under an extra condition imposed on a triangulation, induces cube tiling. We need several technical preparations before getting into the proof. Our proof is a generalization of  the scheme in the two dimensional case given by Schramm (\cite{Sch1}).

\medskip

While a given metric on $G$ is a discrete object in nature, Schramm (\cite[Section 5]{Sch1}) defined a continuous family of metrics which depends on a given curve. 

\begin{Def}
\label{def:familyofmetrics}
Let $\alpha$ be any path in $G$ and let $m$ be any metric on $G$. For $t\geq 0$, we define a one parameter family of metrics on $G$ by
\begin{equation}
m_{\alpha,t}(v)= \left\{
\begin{array}{ll}
  m(v)    & \mbox{for $v\in V\setminus\alpha$}      \\
  m(v)+t & \mbox{for $v\in\alpha$}     
     .\end{array}
\right .
\end{equation}
\end{Def}
In the following, whenever the curve $\alpha$ is specified, we will use the notation $m_t$ instead of  
$m_{\alpha,t}$.
\begin{Lem}
\label{le:dermetric}
For the family of metrics $m_t=m_{\alpha,t}$, $t\in[0,\infty)$ we have
\begin{equation}
\label{eq:derformula}
\left . \frac{d}{dt}(||m_t||_3^{3})\right |_{t=0^{+}}=3\sum_{v\in \alpha}m(v)^2,
\end{equation}
where $\alpha$ is any curve in $G$.
\end{Lem}
\begin{proof}
By definition 
\begin{equation}
\label{eq:m_t}
||m_t||_3^3= \sum_{v\in V}m_t(v)^3=\sum_{v\in V\setminus\alpha}m(v)^3 +\sum_{v\in \alpha}(m(v)^3+3m(v)^2t+ 3m(v)t^2 +t^3).
\end{equation}
The assertion of the lemma follows by subtracting $||m||_3^{3}$ from the right hand-side of the equation above, dividing by $t>0$, and taking the limit as $t\rightarrow 0^{+}$.
\end{proof}

\noindent Let ${\mathcal T}$ be a fixed triangulation of a closed topological cube $Q$, and we let $G=({\mathcal T}^{(0)},{\mathcal T}^{(1)})$ be the corresponding graph.  

For the applications of this paper, we consider the following class of triangulations. 

\begin{Def}
\label{de:intersecting} 
A triangulation of a closed topological cube $Q$ will be said to have the triple intersection property, if the following property of the extremal metric $m_0$ of $(G,B_1,\bar B_1)$ holds. There exist a shortest $m_0$-path joining $B_2$ to $\bar B_2$ and a shortest $m_0$-path joining $B_3$ to $\bar B_3$ which meet all the shortest $m_0$ paths joining $B_1$ to $\bar B_1$. 
\end{Def}

In the figure below the red curves correspond to all the shortest paths joining  $B_1$ to $\bar B_1$, the yellow curve to shortest paths joining $B_2$ to $\bar B_2$,  the green curve to a shortest paths joining $B_3$ to $\bar B_3$; all with respect to the $m_0$ metric.

\begin{figure}[htbp]
\begin{center}
 \scalebox{.40}{ \input{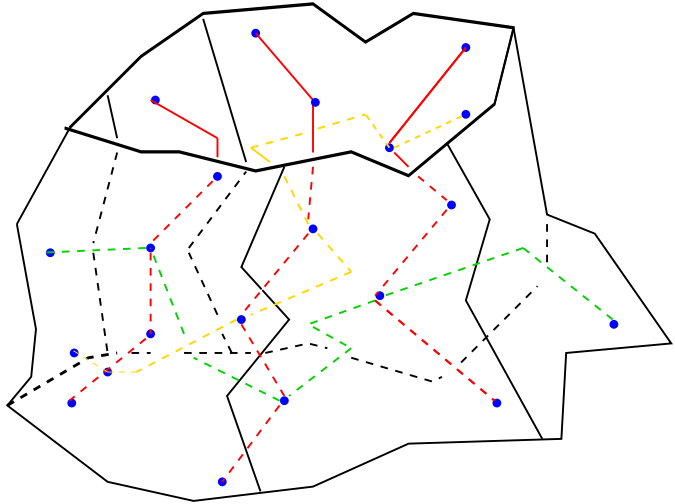}}
 \caption{Shortest paths in a cube.}
\label{figure:quad0}
\end{center}
\end{figure}

Once a triangulation is sufficiently tamed, in the sense described in the following definition, the triple intersection property will hold. It would be interesting to find if  there are other classes of triangulations that have the triple intersection property (see for example Question~\ref{q:q1}). 
We keep the notation of Definition~\ref{de:trianofcube} and make 

\begin{Def}
\label{de:narrowneck}
A triangulation of a closed topological cube $Q$ will be said to have a spine if the following properties hold. There is one and only one path, called the {\it spine} of $G$,  whose interior lies in $B$ and its endpoints lie on $B_1,\bar B_1$, respectively.  In addition, it is required that every 
path whose interior lies in $B$ and joins $B_2$ to $\bar B_2$ or $B_3$ to $\bar B_3$, intersects the {\it spine} of $G$.
\end{Def}

\begin{figure}[htbp]
\begin{center}
 \scalebox{.45}{ \input{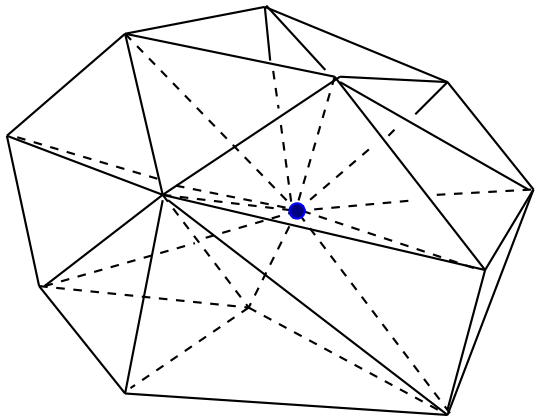}}
 \caption{Part of a triangulation with a spine.}
\label{figure:quad1}
\end{center}
\end{figure}

Suppose that  ${\mathcal T}$ has the triple intersection property, and in addition that $t$ is any nonnegative number which is smaller than the difference between a second shortest $m_0$ path and a shortest $m_0$ path, which join $B_1$ to $\bar B_1$. Consider any $\gamma, \delta$ shortest $m_0$-paths joining $B_2$ to $\bar B_2$ and $B_3$ to $\bar B_3$ as in Definition~\ref{de:intersecting}, respectively. Recall that the metrics $m_{\gamma,t}, m_{\delta,t}$ are obtained by adding $t$ to $m_0(v)$ for each $v\in \gamma$, $v\in \delta$, respectively, and leaving other vertices with their $m_0$ values. Hence, for any $t$ satisfying the condition above, by considering possible shortest paths for $m_{\gamma,t}, m_{\delta,t}$, it follows that 
\begin{equation}
\label{eq:inequalityforpaths}
\min\{l_{m_{\gamma,t}}, l_{m_{\delta,t}}\} \geq l_{m_0} +t,
\end{equation}
and therefore that
\begin{equation}
\label{eq:inequalityforpathsderivative}
\min\{\frac{d}{dt}(l_{m_{\gamma,t}}) |_{t=0^{+}},\frac{d}{dt}(l_{m_{\delta,t}}) |_{t=0^{+}}\}   \geq 1.
\end{equation}

\medskip
Inequality~(\ref{eq:inequalityforpathsderivative}) is essential for the applications of this paper. It will be used in the lemma below (Inequality~(\ref{eq:continue1})) which in turn is essential  in the proof of the main theorem. 
 We continue  with the following lemma which shows that shortest curves measured with respect to the
extremal metric $m_0$ cannot be too short. This will be used in the proof of  Theorem~\ref{th:exttocubes} to show that cubes arising from vertices that belong to boundary components of $Q$ intersect appropriate (extended) boundary components naturally defined by $R$ (for a more precise statement see the proof of Theorem~\ref{th:exttocubes}).

\begin{Lem}
\label{le:shortcurves}
With the notation and hypotheses of Theorem~\ref{th:main} and with $m_0$ being the extremal metric for $(G,B_1,\bar B_1)$ normalized so that  $\mbox{\rm vol}(m_0)=1$, we have
\begin{equation}
\label{eq:notshort}
\min{ \{{ l_{m_{0}}(\alpha), l_{m_{0}}(\beta})}\}\geq \sqrt{h^{-1}},
\end{equation}
where $h=l_{m_0}$, $\alpha$ is any shortest $m_0$ path joining $B_2$ to $\bar B_2$, and $\beta$ is any shortest $m_0$ path joining $B_3$ to $\bar B_3$.
\end{Lem}
\begin{proof}
Since $m_0$ is extremal, we have the following inequality for  $l_{m_t}=l_{m_\alpha,t}$ 
\begin{equation}
\label{eq:firstderivative}
 0\ \geq \ \frac{d}{dt} \left( \frac{(l_{m_t})^3}{||m_t||_{3}^{3}} \right)_{t=0^{+}}=\ \left (\frac{3 l_{m_t}^2  \frac{d}{dt}(l_{m_t})    ||m_t||_3^3-l_{m_t}^{3} \frac{d}{dt}(||m_t||_3^3)}{||m_t||_3^6} \right )_{t=0^{+}}.
\end{equation}

Hence, by applying Lemma~\ref{le:dermetric} and the normalization $\mbox{\rm vol}(m_0)=1$, we must have 
\begin{equation}
\label{eq:continue}
0\geq  \left( l_{m_t}^2 \frac{d}{dt}(l_{m_t}) -   l_{m_t}^3 \sum_{v\in \alpha}m_{0}(v)^2\right )_{t=0^{+}},
\end{equation}
which implies, if $\alpha$ satisfies Inequality~(\ref{eq:inequalityforpathsderivative}), that 
\begin{equation}
\label{eq:continue1}
\sum_{v\in \alpha}m_{0}(v)^2 \geq  l_{m_0}^{-1}=h^{-1}.
\end{equation}

Since  
\begin{equation}
\label{eq:continue2}
l_{m_0}(\alpha)^2=(\sum_{v\in \alpha}m_{0}(v))^2\geq\sum_{v\in \alpha}m_{0}(v)^2,
\end{equation}
and all $m_0$ shortest paths joining $B_2$ to $\bar B_2$ have the same length, the assertion of the lemma follows (an identical argument holds for $\beta$). 

\end{proof}

We do not know if inequality~(\ref{eq:inequalityforpathsderivative}) is necessary for the assertions of Theorem~\ref{th:main} to hold. It is definitely necessary in our proof of the lemma above, as well as in the analogous part of Schramm's proof of his main theorem. The triple intersection property guarantees that this inequality holds.  

\medskip

We now turn into the construction of cube tiling from a normalized extremal metric.

\begin{Thm}
\label{th:exttocubes}
Let ${\mathcal T}=\{V,E,F;B_1,\bar{B_1},B_2,\bar{B_2},B_3,\bar{B_3}\}$ be a triangulation of a topological cube which has the triple intersection property, and let 
$G=(V,E)$ be the $1$-skeleton of ${\mathcal T}$. Let $m$ be the extremal metric for 
$(G,B_1,\bar B_1)$ normalized so that  $\mbox{\rm vol}(m)=1$. Set
\begin{equation}
\label{eq:height} 
h=l_m,\  \mbox{\rm and let }\  R=[0,h]\times[0,\sqrt{h^{-1}}]\times[0,\sqrt{h^{-1}}]. 
\end{equation}
For each $v\in V$ let
\begin{equation}
\label{eq:cubesdimesnions}
C_v=[x(v)-m(v),x(v)]\times[y(v)-m(v),y(v)]\times[z(v)-m(v),z(v)],
\end{equation}
where $x(v)$ $($respectively $y(v), z(v))$  is the least $m$-length of paths from
$\bar B_2$ $($respectively, $B_3,B_1)$ to $v$. Then
${\mathcal C}=\{C_v:v\in V\}$ is a cube tiling of the rectangular parallelepiped $R$  which satisfies the contact  and boundary constraints $(\ref{eq:comb})$ and $(\ref{eq:bound})$.
\end{Thm}
We may now turn to the proof of our main theorem.

\medskip

\noindent{\em Proof of Theorem~\ref{th:main}.} The existence of a cube tiling follows from the existence part in Theorem~\ref{th:existence} and from Theorem~\ref{th:exttocubes}. 
Uniqueness follows from Lemma~\ref{le:cutoex} and the uniqueness part in Theorem~\ref{th:existence}. 
\eop{\ref{th:main}}

We now turn to the 

\medskip

\noindent{\em Proof of Theorem~\ref{th:exttocubes}.}
We start by showing that the combinatorics is preserved in the sense of constraint (\ref{eq:comb}).
Let $(u,v)\in {\mathcal T}^{(1)}$ be given. We claim that 
\begin{equation}
\label{eq:combiispreserved}
x(v)-m(v)\leq x(u)\  \mbox{and}\ x(u)-m(u)\leq x(v).
\end{equation}  
Suppose that $x(v)-m(v)>x(u)$, and let $\alpha_u$ be a shortest $m$ path joining $u$ to $\bar B_2$. 
Then the path $[v,u]\cup \alpha_u$ which joins $v$ to $\bar B_2$ has $m$-length which is equal to  $m(v) + x(u)<x(v)$. This is absurd. Hence, by applying a symmetric argument to prove the second inequality, we have that
\begin{equation}
\label{eq:xcoordiantes}
[x(v)-m(v),x(v)]\cap[x(u)-m(u),x(u)]\neq\emptyset.
\end{equation}  

The argument above goes through for the coordinates $y(v)$ and $z(v)$ (up to replacing $\bar B_2$ with $B_3$ and $B_1$, respectively). Thus, as claimed
\begin{equation}
\label{eq:cubesintersect}
Z_v\cap Z_u\neq\emptyset .
\end{equation} 

We now define several rectangular parallelepipeds in $\RR^3$ where some are degenerate, and the remaining  are infinite and all of which  are naturally associated with $R$.

Let 
$ R_1= \{(x,y,z): \min{\{x,y\}}\geq 0,z = 0\}$, $\hat R_1=\{(x,y,z):\min{\{x,y\}}\geq 0,z\geq h\}$, $R_2 =\{(x,y,z): \min{\{y,z\}}\geq 0, x\geq \sqrt{h^{-1}}\}$, $\hat R_2=\{(0,y,z): \min{\{y,z\}}\geq 0\}$, 
$R_3=\{(x,0,z):  \min{\{x,z\}}\geq 0\}$,  and $\hat R_3=\{(x,y,z): \min{\{x,z\}}\geq 0, y\geq \sqrt{h^{-1}}\}$.

Since the volume of $R$ is equal to one, which by assumption is also equal to $\mbox{\rm vol}(m)$, in order to prove that ${\mathcal C}$ tiles $R$, it suffices to prove the following

\begin{equation}
\label{eq:inclusion}
R\subset \bigcup_{v\in V}C_v, 
\end{equation}
for then it follows that there are no overlaps of positive volume among the cubes and no cube extends beyond $R$.
To this end we first prove    
\begin{Lem}
\label{le:homotopy}
With the notation above we have that $\partial R$ is freely homotopic to a constant in 
\begin{equation}
\label{eq:homotopy}
\bigcup_{v\in V}C_v\cup\left(\RR^3\setminus\mbox{\rm int} (R)\right).
\end{equation}
\end{Lem}

\noindent{\em Proof.}
We begin by constructing a map $f: {\mathcal T}\rightarrow \bigcup_{v\in V}C_v$. For each $v\in V$, choose $f(v)\in C_v$ such that $f(v)\in R_i (\hat R_i)$ for $i=1,2,3$, and whenever $v\in B_i (\bar B_i)$. We observe that this may be done in a consistent way; that is, whenever $v$ is in the intersection of two or three faces among  
$\{B_i,\bar B_i\},$ i=1, 2, 3, then the corresponding intersection among the $\{R_i,\hat R_i\}$ is nonempty. 

We endow ${\mathcal T}$ with a piecewise linear structure by declaring that each $3$-dimensional face $(u,v,w,s)$ of ${\mathcal T}^{(3)}$ is linearly parametrized  by a regular tetrahedron all of its edges have length 1, and that these parametrizations are compatible along faces.  For each $1$-dimensional face $(u,v)\in {\mathcal T}^{(1)}$, let $m_{(u,v)}$ denote  the midpoint of this edge;  for each $2$-dimensional face $(u,v,w)\in {\mathcal T}^{(2)}$,  let $c_{(u,v,w)}$ denote the barycentric center of this face, and for each three dimensional tetrahedron $(u,v,w,s)\in {\mathcal T}^{(3)}$, let $p_{(u,v,w,s)}$ denote its barycentric center.

Choose $f(m_{(u,v)})$ to be some point in $C_v\cap C_u$, choose $f(c_{(u,v,w)})$ to be some point in $C_v\cap C_u\cap C_w$, and choose $f(p_{(u,v,w,s)})$ to be some point in $C_u\cap C_v \cap C_w\cap C_s$. The first choice is possible by applying  (\ref{eq:cubesintersect}), the second and the third choices are possible due to the fact that if three (four) cubes whose edges are parallel to the coordinate axes have the property that the intersection of any two (three) of these cubes is nonempty, then the intersection of the three (four) cubes is nonempty. We also require that $f(m_{(u,v)})\in R_i(\hat R_i)$ if $u,v\in B_i(\bar B_i)$, and that $f(c_{(u,v,w)})\in R_i(\hat R_i)$ if $u,v,w\in B_i(\bar B_i)$. 

Let ${\mathcal T}^{*}$ be the first barycentric subdivision of ${\mathcal T}$, and extend $f$ by requiring it to be affine on each face of  ${\mathcal T}^{*}$. It is clear (by construction) that the extension is well defined. Also, since for each face $(s,m_{s,u},c_{s,u,v},p_{s,u,v,w})$ of  ${\mathcal T}^{*}$ the four points $f(s), f(m_{s,u}), f(c_{s,u,v})$, and  $f(p_{s,u,v,w}))$ lie in $C_s$ which is convex. Hence
\begin{equation}
\label{eq:inclusion1}
f({\mathcal T})\subset \bigcup_{v\in V}C_v.
\end{equation}

Let $v$ be any vertex in ${\mathcal T}^{(0)}$, then it is clear that $\partial {\mathcal T}$ is freely homotopic in ${\mathcal T}$ to $v$. By construction, $f(\partial {\mathcal T})$ is freely homotopic to $\partial R$ in $\RR^3\setminus \mbox{\rm int} (R)$. Note that the last part is justified (in part) due to the assertions of Lemma~\ref{le:shortcurves}. It is here where we are using in an essential way a lower bound for the shortest $m$-curves joining $B_2$ to $\bar B_2$ and $B_3$ to $\bar B_3$. The assertion of the lemma follows immediately by
defining the constant to be $f(v)$, and composing the two homotopies above.
\eop{\ref{le:homotopy}}

We now finish the proof of the theorem by establishing (\ref{eq:inclusion}). We  argue by contradiction. First suppose that there exists a point 
\begin{equation}
x\in \mbox{\rm int}(R)\  \mbox{\rm such that}\ 
x\not\in \bigcup_{v\in V}C_v.
\end{equation} 

Thus, we have the inclusion
\begin{equation}
\bigcup_{v\in V}C_v\cup \left(\RR^3\setminus\mbox{\rm int} (R)\right)\hookrightarrow \RR^3\setminus\{x\}.
\end{equation}

Hence, by the assertion of the previous lemma, $\partial R\simeq \SS^2$ is homotopic to a constant in
$ \RR^3\setminus\{x\}$. This is absurd. To end, one treats the case 
$x \in \partial R$ and $x\not\in \bigcup_{v\in V}C_v$, by arguing that since $\bigcup_{v\in V}C_v$ is a closed set, there exists a point $y\in \mbox{\rm int}(R)$ which is close to $x$ and is not in  $\bigcup_{v\in V}C_v$.

\eop{\ref{th:exttocubes}}

\begin{Rem}
Since at most eight cubes may be tiled in $\RR^3$ with a nonempty intersection, it is feasible that some cubes in the tiling  provided by Theorem~\ref{th:exttocubes} will degenerate to points.
\end{Rem}

\medskip

\begin{Rem}
\label{re:gene}
There are straightforward modifications of our definitions and proofs that allow generalizations of the results to tiling with rectangular parallelepipeds of specified aspect ratios.  In the two dimensional case, one such generalization (tiling by rectangles instead of squares) was observed by Schramm (\cite[Section 8]{Sch1}). 

Let $\omega: V\rightarrow (0,\infty)$ be some assignment of weights to the vertices, and for every metric $m:V\rightarrow [0,\infty)$ on ${\mathcal T}$ define the $\omega m$ length of a path $\alpha=(v_0,v_1,\ldots,v_k)$ as $$l_{\omega m}(\alpha)= \sum_{i=0}^{k}\omega(v)m(v_i),$$ and 
the 
$\omega$-volume  by 
$$ \|m\|_{\omega}^{3}= \sum_{v\in V}\omega(v)m(v)^3 .$$
Define the $\omega$-extremal length of ${\mathcal T}$ to be
$$\lambda(G,B_1,\bar B_1)= \sup_{m}\frac{l_m^3}{\|m\|_{\omega}^{3}},$$ where the supremum is taken over all metrics of positive area. Then, as before $\lambda(G,B_1,\bar B_1)$ is a (discrete) conformal invariant, and a $\omega$-extremal metric always exists, and is unique up to a positive scaling factor.
With this setting, the following holds  (we omit the straightforward details of the proof as well as other possible generalizations).
\begin{Thm}
\label{th;gene}
Let ${\mathcal T}=\{V,E,F;B_1,\bar{B_1},B_2,\bar{B_2},B_3,\bar{B_3}\}$ be a triangulation of a topological cube which has the triple intersection property, and let 
$G=(V,E)$ be the $1$-skeleton of ${\mathcal T}$. Let $\omega: V\rightarrow (0,\infty)$ be some assignment of weights to the vertices.  Let $m$ be the $\omega$-extremal metric for 
$(G,B_1,\bar B_1)$ that satisfies $\|m\|_{\omega}^{3}=1$. Set
\begin{equation}
\label{eq:height1} 
h=l_m,\  \mbox{\rm and let }\  R=[0,h]\times[0,\sqrt{h^{-1}}]\times[0,\sqrt{h^{-1}}]. 
\end{equation}
For each $v\in V$ let
\begin{equation}
\label{eq:cubesdimesnions1}
C_v=[x(v)-\sqrt{\omega(v)}m(v),x(v)]\times[y(v)-\sqrt{\omega(v)}m(v),y(v)]\times[z(v)-m(v),z(v)],
\end{equation}
where $z(v)$ $($respectively, $x(v), y(v))$  is the least $m$-length of a path from $v$ to
$B_1$ (least $\omega m$-length to  $\bar B_2,B_3, respectively)$. Then
${\mathcal C}=\{C_v:v\in V\}$ is a tiling of the rectangular parallelepiped $R$ which satisfies the contact  and boundary constraints $(\ref{eq:comb})$ and $(\ref{eq:bound})$.
\end{Thm}

\end{Rem}

\section{further questions and research directions}
\label{se:section4}

We end this paper by suggesting several future research directions and questions that are motivated in part by the extensive study done in the two dimensional case (see for example \cite{BeSch,De,Ca, CaFlPa,CaFlPa1,Ke}). 

\bigskip

Definition~\ref{de:intersecting} specifies a class of triangulations of a closed topological cube for which the assertions of Theorem~\ref{th:main} holds.
\begin{Que}
\label{q:q1}
Are there larger classes of triangulations which induce a tiling by cubes?  
\end{Que}

\bigskip

Experience shows that the method of extremal length is very useful when two boundary components are chosen  (these are the top base and the bottom base in our work). The passage for $3$-manifolds without boundary invites further investigations.
\begin{Que} 
What is the analogue of Theorem~\ref{th:main} for a ring space domain, and even more generally for a genus $g$ handlebody?
\end{Que}

\bigskip

The works in \cite{Sch1} and in \cite{CaFlPa} contain various algorithms to compute two dimensional extremal length for a triangulation of a quadrilateral.  All of these use the planarity in an essential way. The following seems to be quite natural to pose.
\begin{Que}
Is there an efficient algorithm to compute extremal length for a given ${\mathcal T}$?
\end{Que}

\bigskip

There is a  combinatorial notion of {\it a boundary value} data which may naturally be associated with a cube tiling, that is, the induced square tiling of the faces. We  propose 
\begin{Que} Given a pattern of square tiling of some (perhaps all of) the faces of $R$, does there exist a cube tiling of $R$ that induces this pattern?
\end{Que}

\bigskip
Wood (\cite{Woo,Woo1}) studied how two dimensional discrete extremal length and the associated modulus changes under various effects of combinatorial operations on a triangulated planar annulus, and related questions on triangulated Riemann surfaces.  Without getting into technical definitions, we pose the following.

\begin{Que}
What are the effects of (for example) refinement of a triangulation on the discrete three dimensional extremal length? (We do not have a good understanding of this even in the case discussed in this paper.)
\end{Que}

\bigskip

We close this list of questions by one which is motivated by the classical continuous theory of extremal length. Due to the work of various authors (see the beginning of \S\ref{se:section1}), there are intimate relations between extremal length and harmonic functions. The work in \cite{Her, Her1, Her2} shows that the classical theory does not transform word by word to the discrete setting, tiling by cubes which is induced by harmonic maps is possible, yet more complicated to construct.  
\begin{Que} 
Assume that ${\mathcal T}$ is given (for a topological cube or a handlebody), 
does there exist a tiling by cubes (or by rectangular parallelepipeds) which is  generated by the discrete harmonic function defined on $V$ and which satisfies suitable combinatorial boundary conditions (such as Dirichlet or Dirichlet-Neumann)?
\end{Que}


\begin{thebibliography}{99}


\bibitem{Ah} L.~V.~Ahlfors, \emph{Conformal invariants-Topics in Geometric Function Theory}, McGraw-Hill Book Company, 1973.

\bibitem{An1} E.M.~Andreev, \emph{On convex polyhedra in Loba\u{c}evski\u{i} space}, Mathematicheskii Sbornik (N.S.) \textbf{81} (\textbf{123}) (1970), 445--478 (Russian); Mathematics of the USSR-Sbornik \textbf{10} (1970), 413--440 (English).

\bibitem{An2} E.M.~Andreev, \emph{On convex polyhedra of finite volume in Loba\u{c}evski\u{i} space},  
Mathematicheskii Sbornik (N.S.) \textbf{83} (\textbf{125}) (1970), 256--260 (Russian); Mathematics of the USSR-Sbornik \textbf{10} (1970), 255--259 (English).

\bibitem{BeSch} I.~Benjamini and O.~Schramm, \emph{Random walks and harmonic functions on infinite planar graphs using square tilings}, Ann. Probab. \textbf{24} (1996), 1219--1238.

\bibitem{BeSch1} I.~Benjamini and O.~Schramm, \emph{Lack of Sphere Packing of Graphs via Non-Linear Potential Theory}, preprint, arXiv:0910.3071v2.

\bibitem{BeCu} I.~Benjamini and N.~Curien, \emph{On limits of Graphs Sphere Packed in Euclidean Space and Applications}, preprint, arXiv:0907.2609v4.

\bibitem{BSST} R.L~Brooks, C.A.~Smith, A.B.~Stone and W.T.~Tutte, \emph{The dissection of squares into squares}, Duke Math. J. \textbf{7} (1940), 312--340.

\bibitem{Ca} J. W.~Cannon, \emph{The combinatorial Riemann mapping theorem},  Acta Math. \textbf {173} (1994), 155--234. 


\bibitem{CaFlPa} J.W.~Cannon, W.J.~Floyd and W.R.~ Parry,  \emph{Squaring rectangles: the finite Riemann mapping theorem},  Contemporary Mathematics, vol. 169, Amer. Math. Soc., Providence, 1994, 133--212.

\bibitem{CaFlPa1} J.W.~Cannon, W.J.~Floyd and W.R.~ Parry, \emph{Squaring rectangles for dumbbells}, Conform. Geom. Dyn. \textbf{12} (2008), 109--132

\bibitem{De} M.~Dehn, \emph{Zerlegung ovn Rechtecke in Rechtecken}, Mathematische Annalen, \textbf {57}, (1903), 144-167.

\bibitem{Ke} R.~Kenyon, \emph{Tilings and discrete Dirichlet problems}, Israel J. Math. \textbf{105} (1998), 61--84.

\bibitem{Ko} P.~Koebe, \emph{Kontaktprobleme der konforrmen Abbidung}, S\"{a}chs. Akad. Wiss. Leipzig, Math. -Phys. Klasse \textbf{88} (1936), 39--53. 

\bibitem{Lo} C. L\"oewner, \emph{On the conformal capacity in space}, Jour. Math. Mech. {\textbf 8} (1959), 411--414.

\bibitem{RoSu} B.~Rodin and D.~Sullivan, \emph{The convergence of circle packings to the Riemann  mapping}, Jour. Differ. Geometry \textbf{26} (1987), 349--360.

\bibitem{Ge} F.W.~Gehring, \emph{Extremal length definitions for the conformal capacity of rings in space}, Michigan Math. Jour. \textbf{9} (1962), 137--150.


\bibitem{Hei} J.~Heinonen, \emph{Lectures on analysis on metric spaces}, Universitext,  Springer-Verlag, New York, 2001.

\bibitem{Her} S.~Hersonsky, \emph{Energy and length in a topological planar quadrilateral},  
European J. Combin. \textbf{29} (2008), no. 1, 208--217.

\bibitem{Her1} S.~Hersonsky, \emph{Boundary Value Problems on Planar Graphs and Flat Surfaces with Integer Cone singularities I; The Dirichlet problem}, Crelle, accepted  for publication.

\bibitem{Her2} S.~Hersonsky, \emph{Boundary Value Problems on Planar Graphs and Flat Surfaces with Integer Cone singularities II; Dirichlet-Neumann problem},  Differential geometry and its applications. accepted for publication.

\bibitem{Sa1} B.V.~\u{S}abat, \emph{The modulus method in space}, Soviet Math. {\textbf 1} (1960), 165--168.

\bibitem{Sa2} B.V.~\u{S}abat, \emph{On the theory of quasiconformal mappings in space}, Soviet Math. {\textbf 1} (1960), 730--733.


\bibitem{Sch1} O.~Schramm, \emph{Square tilings with prescribed combinatorics}, Israel Jour. of Math. \textbf{84} (1993), 97--118.


\bibitem{Sch2} O.~Schramm, \emph{Conformal uniformization and packings}, Israel Jour. of Math. \textbf{93} (1996), 399--428.

\bibitem{Th1} W.P.~Thurston, \emph{The finite Riemann mapping theorem}, invited address, International Symposium in Celebration of the Proof of the Bieberbach Conjecture, Purdue University, 1985.

\bibitem{Th2} W.P.~Thurston, \emph{The Geometry and Topology of $3$-manifolds}, Princeton University Notes, Princeton, New Jersey, 1982.

\bibitem{Va1} J.~V\"{a}is\"{a}l\"{a}, \emph{On quasiconformal mappings in space}, Ann. Acad. Sci. Fenn. Ser. A. I. \textbf{298} (1961), 1--36.


\bibitem{Va2} J.~V\"{a}is\"{a}l\"{a}, \emph{On quasiconformal mappings of a ball}, Ann. Acad. Sci. Fenn. Ser. A. I. \textbf{304} (1961), 1--7.

\bibitem{Woo} William E.~Wood, \emph{Combinatorial modulus and type of graphs}, Topology and its Applications, \textbf{156} (2009), 2747--2761.

\bibitem{Woo1} William E.~Wood, \emph{Bounded outdegree and extremal length on discrete Riemann surfaces}, Conform. Geom. Dyn. \textbf{14} (2010), 194--201.



\end{thebibliography}
\end{document}